\pdfoutput=1
\documentclass{article}

\usepackage{arxiv}
\usepackage{amsmath}
\usepackage[utf8]{inputenc} 
\usepackage[T1]{fontenc}    
\usepackage{hyperref}       
\usepackage{cleveref}
\usepackage{url}            
\usepackage{booktabs}       
\usepackage{amsfonts} 
\usepackage{float}
\usepackage{nicefrac}       
\usepackage{microtype}      
\usepackage{lipsum}
\usepackage{graphicx}
\usepackage{subcaption}
\usepackage{natbib}
\graphicspath{ {./images/} }

\title{Hub Network Design Problem with Capacity, Congestion and Heterogeneous Economies of Scale}

\author{
 Xiaotong Liu \\
  Tsinghua-Berkeley Shenzhen Institute\\
  Tsinghua Shenzhen International Graduate School\\ Tsinghua University\\
  \texttt{xiaotong18@mails.tsinghua.edu.cn} \\
}

\begin{document}
\maketitle
\begin{abstract}
We propose a joint model that links the strategic level location and capacity decisions with the operational level routing and hub assignment decisions to solve hub network design problem with congestion and heterogeneous economics of scale. We also develop a novel flow-based mixed-integer second-order cone programming (MISOCP) formulation. We perform numerical experiments on a real-world data set to validate the efficiency of solving the MISOCP reformulation. The numerical studies yield observations can be used as guidelines in the design of transportation network for a logistics company.
\end{abstract}


\section{Introduction}\label{sec4-introduction}
Hub-and-spoke topologies are widely applicable in a variety of networks such as airline transportation \citep{shen2021reliable}, postal delivery \citep{ernst1999solution}, express delivery service \citep{yaman2012release,wu2023service}, telecommunication \citep{klincewicz1998hub}, and brain connectivity networks \citep{khaniyev2020spatial}. Hubs help reduce the cost of establishing a network connecting many origins and destinations, and also consolidate flows to exploit economies of scale. A hub location problem is a network design problem consisting generally of two main decisions: to locate a set of hubs and to allocate the demand nodes to these hubs \citep{campbell2012twenty}.

To capture the consolidation-deconsolidation trade-off, the classical hub location problem assume the unit transportation cost on the intermediate link between two hubs is discounted by a constant factor independent of the actual flow. In practice, however, as the flow increases, the unit transportation cost decreases due to sharing the fixed costs over more units of goods and the potential use of cost efficient vehicles (e.g. aircraft and trucks). As such, many researchers began to consider flow-dependent or heterogeneous economies of scale on interhub flow. \cite{o1998hub} assumed the cost function to be a piecewise-linear function of the flow, with each flow segment having different fixed and variable cost. 

Considering heterogeneous economies of scale creates a force pushing for consolidation. However, over-consolidation can create congestion at hub nodes, which may lead to high overall cost and poor service quality \citep{alumur2018modeling}. In particular. when the traffic flow through a hub approaches its capacity, the operational cost incurred due to congestion, increases steeply. Therefore, nonlinear modeling of the congestion cost yields more realistic results. Several studies in the literature have explicitly considered congestion costs and capacity acquisition decisions in the hub location problem \citep{elhedhli2010lagrangean,alumur2018modeling,bayram2023hub}. All of the above works, however, only adopt the simplistic flow-independent method described above. This fails to capture the interplay between consolidation and congestion. To the best of our knowledge, this paper is the first one to propose a joint model that links the strategic level location and capacity decisions with the operational level routing and hub assignment decisions to solve hub network design problem with congestion and heterogeneous economics of scale.  

 We introduce the hub network design problem with capacity, congestion and heterogeneous  economics of scale (HNDCH) to cover a broad range of strategic to operational level decisions. HNDCH aims to find the optimal network design that minimizes the total setup, capacity acquisition, congestion and routing (transportation) cost. The introduced problem has the following key characteristics. (1) Hubs are capacitated, that is, the total flow through a hub location is restricted. The network management incurs a congestion cost that depends on how much of the available capacity is used at a hub. (2) Heterogeneous economies of scale on interhub flow. Instead of a constant multiplier, interhub flow costs were integrated into the cost through a piecewise-linear concave function.

 To model this challenging problem, we propose a novel flow-based mixed-integer second-order cone programming(MISOCP) formulation. We aim to find the optimal solution to the HNDCH and answer the following research question. 1)How does the optimal network topology change when we consider the congestion cost and the heterogeneous economics of scale? To this end, extensive computational results are conducted on a large real-world data set provided by a logistics company, which contains 688 nodes (30 of which are potential hub locations).

The rest of the paper is organized as follows. Section \ref{sec4-preliminaries} describes the problem formulation. Section \ref{sec4-MISOCP} details the flow-based MISOCP reformulation. Section \ref{sec4-computer} discusses the computational results and presents the managerial insights, and Section \ref{sec4-conclude} concludes.

 
\section{Preliminaries and Model}\label{sec4-preliminaries}
This section first introduces some preliminaries needed throughout this chapter and then presents mathematical formulations for the HNDCH problem.

\subsection{Network Topology}\label{sec4-top}
Consider a directed graph $G = (N, A)$ where $N = \{1, \dots, n\}$ is the set of nodes and $A$ is the set of arcs in the network. From this point on, we use indices $i, j, k, m \in N$ for nodes and drop the set notation for convenience. Let $w_{ij}$ be the amount of flow to be transported from node $i$ to node $j$, and let $d_{ij}$ be the distance from node $i$ to node $j$. We define $O_i = \sum_j w_{ij}$ and $D_i = \sum_j w_{ji}$ as the total outgoing flow from node $i$ and the total incoming flow to node $i$, respectively. 

We define a path $p = (i, k, m, j)$ from an origin node $i$ to a destination node $j$ passing through hubs $k$ and $m$ in that order. We assume that at most one hub arc will be used in each path. In other words, for each path $p = (i, k, m, j)$ with $k \neq m$, the traffic $w_{ij}$ is sent on arc $(i, k)$, is then routed through the interhub connection $(k, m)$, and is finally delivered on arc $(m, j)$. If $k = m$, the traffic $w_{ij}$ flows on arc $(i, k)$ and then on arc $(k, j)$. 

\begin{figure}
  \centering
  \subcaptionbox{path $i-k-m-j$}
    {\includegraphics[width=0.42\textwidth]{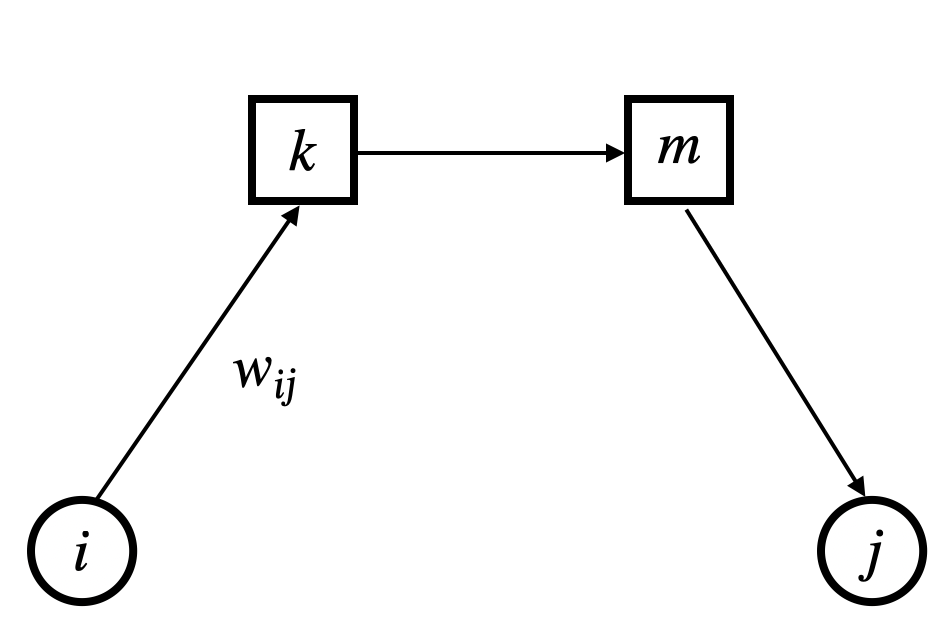}}
  \subcaptionbox{path $i-k-j$}
    {\includegraphics[width=0.38\textwidth]{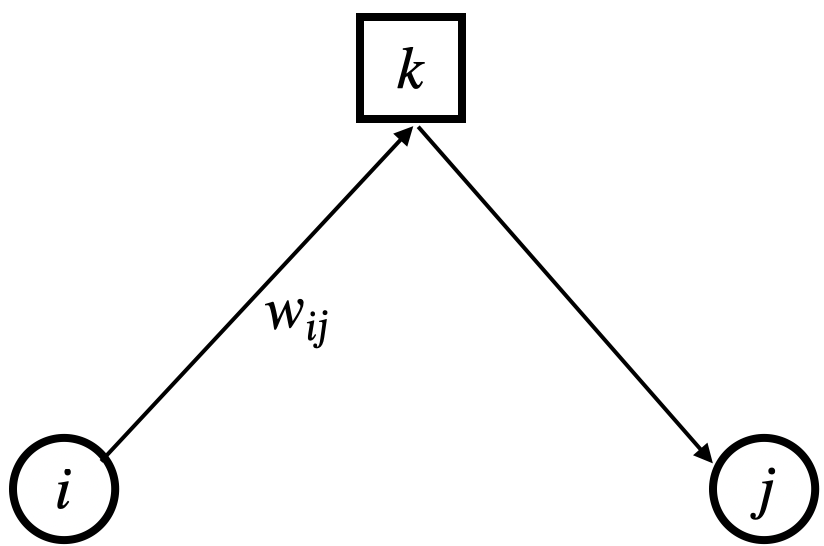}}
  \caption{Illustration of Paths}
  \label{fig-hub-path}
\end{figure}

\subsection{Heterogeneous Economies of Scale}\label{sec4-scale}
There are three types of cost associated with each path $p$: the collection cost $\chi_{ik}$, the transfer cost $F_{km}$, and the distribution cost $\delta_{m j}$. We assume that $\chi_{ik}$, $F_{km}$, and $\delta_{m j}$ are concave functions of the total flow on arcs $(i, k)$, $(k, m)$, and $(m, j)$, respectively. Because of the single allocation assumption, the values of the collection and distribution functions $\chi_{ik}$ and $\delta_{m j}$ depend only on the flows $O_i$ and $D_i$ and hence can be determined a priori. For notational convenience, let $c_{ik} = \chi_{ik} + \delta_{ki}$ be the total collection and distribution cost between node $i$ and hub $k$. Let $v_{km}$ be the amount of flow routed on the interhub connection $(k, m)$ with $k \neq m$ and $v_{kk} = 0$ for all $k$ because the cost associated with the flow between two spoke nodes that are assigned to the same hub node are captured by the $c_{ik}$ parameters. The term $F_{km}$ is then a concave function over the finite interval $[a_{km}, b_{km}]$ where $a_{km}$ and $b_{km}$ are the lower and upper bounds on the flow variable $v_{km}$.

\begin{figure}
    \centering
    \includegraphics[width=0.6\textwidth]{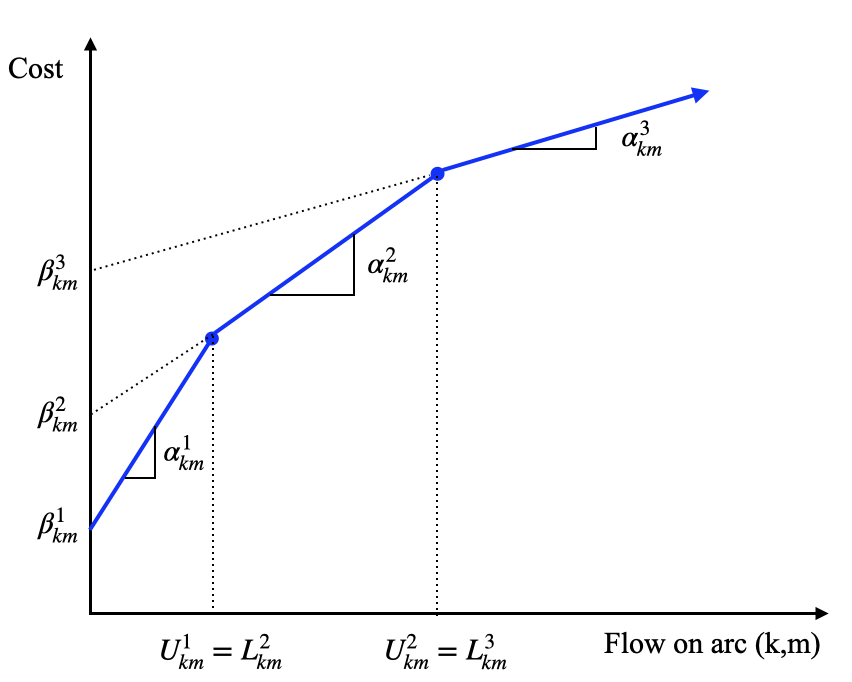}
    \caption{Piecewise-Linear Concave Cost Function on arc $(k,m)$}
    \label{fig:Hub-scale}
\end{figure}

Taking into account the above mentioned assumptions and requirements on the structure of origin-destination
paths, we define the cost function $F_{km}: s \in S \rightarrow \mathbb{R}$ with $I_s = [L^s_{km}, U^s_{km}]$ as:
\begin{equation}
    F_{km}(v_{km}) = \min_{s \in S} (\beta^s_{km} + \alpha^s_{km} v_{km}) \quad \forall k,m,
\end{equation}
where $\beta^s_{km}$ and $\alpha^s_{km}$ are the fixed and variable costs of operating in segment $s$ on arc $(k, m)$, respectively. The coefficients $\beta^s_{km}$ are increasing with $s$ and the coefficients $\alpha^s_{km}$ are decreasing with $s$.

\subsection{Capacity and Congestion}\label{sec4-capacity}
A hub can be opened in different sizes (capacity levels) listed in the set $L$. The fixed and variable cost of opening a size $\ell \in L$ hub at location $k\in N$ is given by $f_k^\ell$. The capacity associated with a hub size $\ell$ is denoted by $q^\ell$ and the maximum possible capacity for a hub (the capacity of the largest size hub) is shown as $Q$.

\begin{figure}
    \centering
    \includegraphics[width=0.6\textwidth]{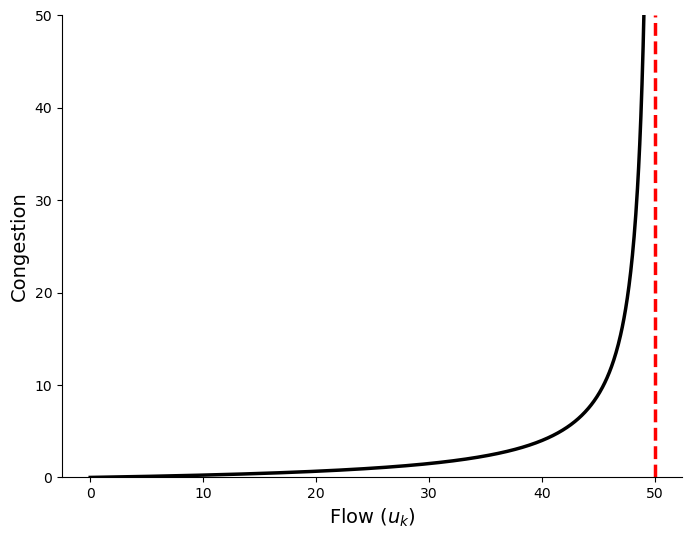}
    \caption{Illustration of Kleinrock Congestion Function at Hub $k$}
    \label{fig:Hub-congestion}
\end{figure}

To reflect the congestion cost, a Kleinrock function is used that models each hub as an $M/M/1$ queue in the steady-state condition \citep{kleinrock2007communication, elhedhli2010lagrangean}. Two 
input elements are critical for queuing analysis at a hub: mean arrival flow rate and mean service rate (i.e.,
vehicles per hour). The service capacity for transshipment flows at hub, denoted as $C_k$, is comparable to the mean service rate, whereas $u_k$, the total transshipment flow at hub $k$, is analogous to the mean arrival flow rate. Hence, it is appropriate to formulate the congestion rate at a hub $k$ as 
\begin{equation}\label{eq:congestion}
    \frac{u_k}{C_k-u_k} \quad \forall k.
\end{equation}
As shown in Figure \ref{fig:Hub-congestion}, the closer the total flow to the capacity (red line), the more severe the congestion is (black line). 

\subsection{Model}\label{sec4-model}
The HNDCH links strategic level hub location and capacity decisions with operational level hub assignment and routing decisions to achieve an optimal hub network design, that is, number, location, and capacity selection, congestion, and transportation cost is minimized.  To formulate HNDCH, we introduce the variables:
\begin{itemize}
    \item $x_{ik} = 1$ if node $i$ is allocated to a hub located at $k$, $0$ otherwise.
    \item For every node $k$, $x_{kk} = 1$ indicates whether $k$ is a hub ($x_{kk} = 1$) or not ($x_{kk} = 0$).
    \item We define a binary variable $t_k^\ell$ for each $k$ and $\ell \in L$ as follows:
    \[
    t_k^\ell = \begin{cases}
    1 & \text{if a hub $k$ of size $\ell \in L$ is open,}\\
    0 & \text{otherwise.}
    \end{cases}
    \]
    \item We also define a binary variable $z^s_{km}$ for each $k, m$ and $s \in S$ as follows:
    \[
    z^s_{km} = \begin{cases}
    1 & \text{if the amount of flow between hubs $k$ and $m$ falls in the interval $I_s$,} \\
    0 & \text{otherwise.}
    \end{cases}
    \]
\end{itemize}
For convenience, we refer to $x$ as hub location variables, to $t$ as capacity selection variables and to $z$ as flow segment variables. We use an $s$ index for the flow segment variables. Consequently, $v_{km}$ and $F_{km}(v_{km})$ can be expressed as
\begin{align}
    v_{km} &= \sum_{i, j} w_{ij} x_{ik} x_{jm} \quad \forall k, m, k \neq m, \\
    F_{km}(v_{km}) &= \sum_{s} z^s_{km} (\beta^s_{km} + \alpha^s_{km} v_{km}),
\end{align}
where the first equation follows from the single allocation assumption, and the second equation follows from the definition of the $z$ variable by which at most a single flow segment variable can be nonzero on any hub-hub connection.

Moreover, $u_k$ and $C_k$ in \eqref{eq:congestion} can be expressed as
\begin{align}
    u_k  &= \sum_{i} \sum_{j} w_{ij}x_{ik} = \sum_{i} O_i x_{ik} \quad \forall k , \\ 
    C_k &= \sum_{\ell}q^\ell t_k^\ell, \quad \forall  k,
\end{align}
and the congestion cost function for hub $k$ is given by 
\begin{equation*}
g_k \frac{\sum_{i} O_ix_{ik}}{\sum_{\ell}q^\ell t_k^\ell- \sum_{i} O_i x_{ik}} \quad \forall  k,
\end{equation*}
where $g_k$ is a scaling factor used to calculate the congestion of hub $k$. The HNDCH can then be formulated as the following mixed integer nonlinear program (MINLP):
\begin{align}
\min \quad & \sum_{k}\sum_{\ell} f_k^\ell t_k^\ell + \sum_{k}g_k \frac{\sum_{i} O_ix_{ik}}{\sum_{\ell}q^\ell t_k^\ell- \sum_{i} O_i x_{ik}}+\sum_{i,k} c_{ik} x_{ik} \label{eq:object1}\\
&+\sum_{k, m, s} d_{km} \beta^s_{km} z^s_{km} + \sum_{i, j, k, m, s} d_{km} \alpha^s_{km} z^s_{km} w_{ij} x_{ik} x_{jm} \nonumber\\
\text{s.t.} \quad & \sum_{k} x_{ik} = 1 & \forall i, & \label{eq:allocation}\\
& x_{ik} \leq x_{kk} & \forall i, k, & \label{eq:hub_existence} \\
& \sum_{i} O_{i} x_{ik} \leq \sum_{\ell}q^\ell t_k^\ell & \forall k,  & \label{eq:capacity} \\
& \sum_{i,j} w_{ij} x_{ik} x_{jm} \geq \sum_{s}L^s_{km} z^s_{km} & \forall k, m, k \neq m,  & \label{eq:flow_lower_bound} \\
& \sum_{i,j} w_{ij} x_{ik} x_{jm} \leq \sum_{s}U^s_{km} z^s_{km} & \forall k, m, k \neq m,  & \label{eq:flow_upper_bound} \\
& \sum_{s} z^s_{km} \geq x_{kk} + x_{mm} - 1 & \forall k, m, k \neq m, & \label{eq:segment_activation} \\
& \sum_\ell t_k^\ell = x_{kk} & \forall k, & \label{eq:capacity_activation} \\
& x_{ik} \in \{0, 1\} & \forall i, k, & \label{eq:binary_allocation} \\
& z^s_{km} \in \{0, 1\} &\forall k, m, s, & \label{eq:binary_flow_segments}\\
& t_k^\ell \in \{0, 1\} & \forall k, \ell. & \label{eq:binary_capacity_selection}
\end{align}
The objective is to minimize the total cost that includes four components: the cost of opening new hubs with a certain capacity, the congestion costs, the cost of assigning the spoke nodes to the hub nodes, and the cost of transporting goods on hub-hub lines. 
Constraints \eqref{eq:allocation} guarantee that each node is assigned to exactly one hub, whereas
constraints \eqref{eq:hub_existence} impose that they can only be assigned to open hubs. Constraints \eqref{eq:capacity} restrict the total flow incoming to hubs. Constraints \eqref{eq:flow_lower_bound} and \eqref{eq:flow_upper_bound} ensure that the flow on each interhub arc $(k, m)$ lies within the interval $I_s$ if $z^s_{km} = 1$. Constraints \eqref{eq:segment_activation} force the activation of one segment $s$ for each arc $(k, m) \in A$ if both nodes $k$ and $m$ are selected as hub nodes. Constraints \eqref{eq:capacity_activation} state that when a hub is established it should be set up with exactly one capacity level. Constraints \eqref{eq:binary_allocation}- \eqref{eq:binary_capacity_selection} define variable domains.

\section{Model Reformulation}\label{sec4-MISOCP}
The objective function of HNDCH is to minimize a nonconvex function over a nonconvex set, which makes the problem intractable for any of the available solvers. There are two sources of nonlinearity in our model, including the quadratic terms $x_{ik}x_{jm}$ in constraints \eqref{eq:flow_lower_bound} and \eqref{eq:flow_upper_bound}, the congestion term and the cubic terms $z^s_{km} x_{ik} x_{jm}$ in the objective function.

In this section, we first introduce a flow-based reformulation to lineraize the quadratic and the cubic terms. Then we discuss the SOCP transformation to abtain the MISOCP formulation for the problem with a linear objective function and SOCP constraints.

\subsection{Flow-based Reformulation}
As discussed in \cite{campbell2005hub}, the flow-based formulations use continuous variables to compute the amount of flow routed on a particular arc originated at a given node. Under the assumption of single assignments, we only need to use one set of flow variables associated with the hub arcs.  
Following \cite{rostami2022single}, we define a new variable $y_{ikm}^{ s}$ for each $i, k, m, s$ as the total amount of flow originating at node $i$ and routed via hubs located at nodes $k$ and $m$ using segment $s$. 

Updating the constraints \eqref{eq:flow_lower_bound} and \eqref{eq:flow_upper_bound} yields the inequalities
\begin{align}
    & \sum_i y_{i k m}^s \leq U_{k m}^s z_{k m}^s \quad \forall s, k, m, k \neq m, \label{eq:upper_bound2}\\
    & \sum_i y_{i k m}^s \geq L_{k m}^s z_{k m}^s \quad \forall s, k, m, k \neq m.\label{eq:lower_bound2}
\end{align} 
Because of the concave cost structure, Constraints \eqref{eq:lower_bound2} are not required. Then, HNDCH can be reformulated as the following flow-based MINLP:
\begin{align}
\text{MINLP-flow:} \quad \min \quad & \sum_{k}\sum_{\ell} f_k^\ell t_k^\ell + \sum_{k}g_k \frac{\sum_{i,m,s} y_{ikm}^s}{\sum_{\ell}q^\ell t_k^\ell- \sum_{i,m,s} y_{ikm}^s}+\sum_{i,k} c_{ik} x_{ik} \label{eq:MINLP-flow-objective}\\
&+\sum_{k, m, s} d_{km} \beta^s_{km} z^s_{km} + \sum_{i, k, m, s} d_{km} \alpha^s_{km} y_{ikm}^s \nonumber\\
\text{s.t.} \quad & \eqref{eq:allocation}-\eqref{eq:capacity}, \eqref{eq:segment_activation}-
\eqref{eq:upper_bound2},\nonumber 
\\
&\sum_{s }\sum_{m } Y^s_{i k m}=O_i x_{i k} \quad\forall i,k,\label{eq:flow_conservation1}
\\
&\sum_{s }\sum_{m} Y^s_{i m k}=\sum_{j} w_{i j} x_{j k} \quad\forall i,k,\label{eq:flow_conservation}
\\
& y_{ikm}^s\geq 0, \quad \forall i,k,m,k\neq m, s. \label{eq:flow_variable}
\end{align}
Constraints \eqref{eq:flow_conservation1} and \eqref{eq:flow_conservation} are the flow conservation constraints for each O/D node $i$ at each (potential) hub node $k$, where the supply and demand at each node depends on the allocation decision. In this way, the quadratic and cubic term in the Objective \eqref{eq:object1} is  linearized. However, the congestion term in Objective \eqref{eq:MINLP-flow-objective} is nonconvex function. Hence, we make the following modifications to the model.

We define $u_k^\ell$ as the total flow through a hub at location $k \in N$ with capacity $\ell \in L$. Thus 
\begin{equation} \label{def:flow_location_capacity}
   \sum_{\ell} u_k^\ell = \sum_{i,m,s} y_{ikm}^s \quad \forall  k.
\end{equation}
We replace the Objective \eqref{eq:MINLP-flow-objective} with the following convex function:
\begin{align}
    \min \quad & \sum_{k}\sum_{\ell} f_k^\ell t_k^\ell + \sum_{k}\sum_{\ell}g_k \frac{u_k^\ell}{q^\ell - u_k^\ell}+\sum_{i,k} c_{ik} x_{ik} \nonumber\\ &+\sum_{k, m, s} d_{km} \beta^s_{km} z^s_{km} + \sum_{i, k, m, s} d_{km} \alpha^s_{km} y_{ikm}^s.\nonumber
\end{align}
We update the Constraint \eqref{eq:capacity} with \eqref{def:flow_location_capacity}, which we define as follows:
\begin{align}
    &u_k^\ell \leq q^\ell t_k^\ell \quad \forall k,\ell, \label{eq:capacity3}\\
    &u_k^\ell \geq 0 \quad \forall k,\ell. \label{eq:capacity_variable}
\end{align}
With the state updates, the resulting problem is the convex problem. 
\subsection{Mixed-Integer Second-Order Cone Programming Formulation}
Here, we reformulate the HNDCH problem as a MISOCP where the nonlinearity is transferred from the objective function to the constraint set in the form of second order quadratic constraints. The advantage of the MISOCP formulation is that it can be solved directly using  standard  optimization  software  packages  such  as CPLEX or Mosek. To achieve this, we define an auxiliary variable $r_k^\ell, k,\ell$ as follows:
\begin{align}
        &r_k^\ell \geq 0   \quad \forall k, \ell, \label{eq:auxiliary variable}\\
        &r_k^\ell \geq \frac{u_k^\ell}{q^\ell-u_k^\ell}  \quad \forall k, \ell. \label{eq:auxiliary variable2}
\end{align}
We transform \eqref{eq:auxiliary variable2} into a second-order cone constraint by multiplying both sides of it by $q^\ell$ and adding $(u_k^\ell)^2$ to both sides, which yields
\begin{equation}\label{eq:conic}
    (u_k^{\ell})^2 \leq (q^{\ell}r_k^\ell - u_k^{\ell})(q^{\ell} - u_k^{\ell}) \quad \forall k, \ell.
\end{equation}
Constraint \eqref{eq:conic} is a \emph{hyperbolic inequality} of the form $\zeta \leq \xi_{1}\xi_{2}$, where $\zeta, \xi_{1}, \xi_{2} \geq 0$. The constraint $\zeta^2 \leq \xi_{1}\xi_{2}$ can be transformed into the quadratic form $\| (2\zeta, \xi_{1} - \xi_{2}) \| \leq \xi_{1} + \xi_{2}$, where $\| \cdot \|$ is the Euclidean norm \citep{alizadeh2003second}. Following \cite{gunluk2008perspective}, we can represent Constraint \eqref{eq:conic} as the following second-order cone constraint:
\begin{equation}\label{eq:conic2}
\left\|\left(\begin{array}{c}
2 u_k^\ell \\
q^{\ell}r_k^\ell - q^{\ell}
\end{array}\right)\right\| \leq q^{\ell}r_k^\ell + q^{\ell} - 2u_k^\ell \quad \forall k, \ell. 
\end{equation}
Using the previous transformations, we can formulate the HNDCH as the following flow-based MISOCP:
\begin{align}
\text{MISOCP-flow:} \quad \min \quad & \sum_{k}\sum_{\ell} f_k^\ell t_k^\ell + \sum_{k}\sum_{\ell}g_k r_k^\ell+\sum_{i,k} c_{ik} x_{ik} \label{eq:MISOCP-flow-objective}\\ &+\sum_{k, m, s} d_{km} \beta^s_{km} z^s_{km} + \sum_{i, k, m, s} d_{km} \alpha^s_{km} y_{ikm}^s\nonumber\\
\text{s.t.} \quad & \eqref{eq:allocation}-\eqref{eq:hub_existence}, \eqref{eq:segment_activation}-\eqref{eq:upper_bound2},\\
&\eqref{eq:flow_conservation}-\eqref{eq:auxiliary variable}, \eqref{eq:conic2}\nonumber .
\end{align}
This MISOCP-flow model contains a polynomial number of variables and linear constraints plus second-order cone constraint. Hence, it can be solved directly through a commercial solver like Mosek.

\section{Computational Results}\label{sec4-computer}
In  this  section  we  present  our  computational  results  on solving  the  corresponding MISCOP formulations of the joint HNDCH problems discussed in the previous sections. We perform our computational test on Colab with 50.1 GB RAM. The algorithm is coded in Python using Mosek as the mathematical solver. The data is provided by a logistics company, which includes 2,672 routes connecting 642 sources to 16 destinations. To develop less-than-truckload (LTL) transportation, the company is considering leasing warehouses for freight consolidation and transshipment. The parameters of the warehouses are as follows:
\begin{table}[H]
    \centering
    \caption{Hub parameters}
    \begin{tabular}{c|c|c|c}
    \toprule
       Parameter  &  Size 1   &  Size 2 &  Size 3 \\
    \midrule
        $q^\ell$ &     10000 &  15000 &  20000\\
        $f^\ell$ &     12.5&  12.5 &  10\\
        \bottomrule
    \end{tabular}
    \label{tab:hub}
\end{table}

Furthermore, the company possesses a fleet of vehicles comprising two distinct sizes, which can represent the heterogeneous economics of scales. Considering the vehicle specifications, we establish the transportation costs between hubs as below:
\begin{table}[H]
    \centering
    \caption{Segment parameters}
    \begin{tabular}{c|c|c}
    \toprule
       Parameter  &  Size 1   &  Size 2    \\
    \midrule
        $\beta^s$ &     500 &  800 \\
        $\alpha^s$ &     0.108&  0.056 \\
        $U^s $ &     72&  126 \\
    \bottomrule
    \end{tabular}
    \label{tab:segment}
\end{table}

\begin{figure}
  \centering
  \subcaptionbox{heterogeneous economics of scales }
    {\includegraphics[width=0.45\textwidth]{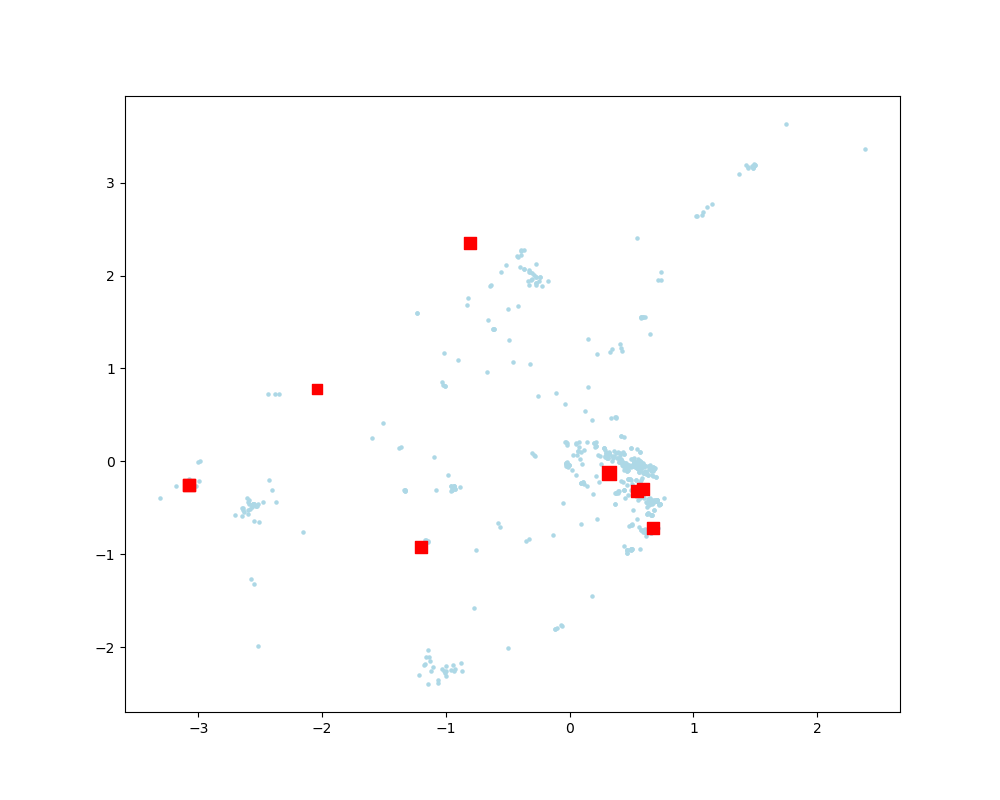}}
  \subcaptionbox{homogeneous economics of scales}
    {\includegraphics[width=0.45\textwidth]{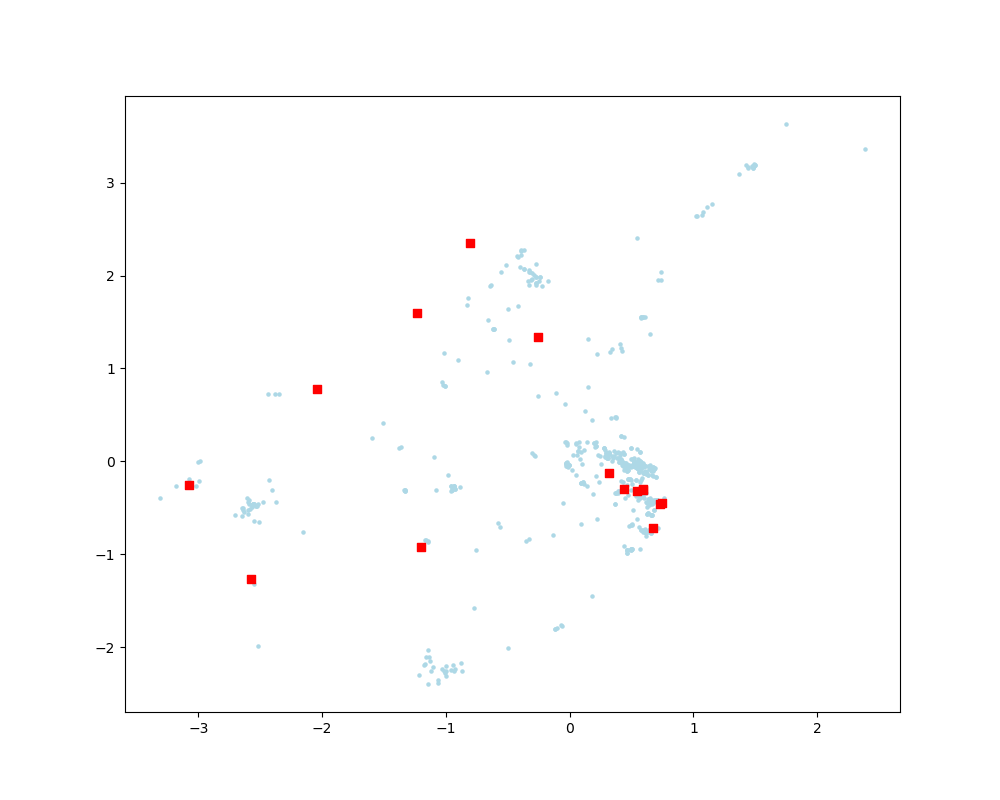}}
  \caption{Hub locations} \label{fig-hub}
\end{figure}

We investigate the hub network topology differences with and without heterogeneous economics of scales consideration. Figure \ref{fig-hub} show that when we consider the heterogeneous economics of scales, only seven hubs with big capacity is able to serve the whole network. However, when we consider homogeneous economics of scales, the resulting design opens a large number of hubs with smaller capacity.

\section{Conclusion}\label{sec4-conclude}
This study propose a joint model that links the strategic level location and capacity decisions with the operational level routing and hub assignment decisions to solve hub network design problem with congestion and heterogeneous economics of scale. We also develop a novel flow-based mixed-integer second-order cone programming (MISOCP) formulation. We perform extensive numerical experiments on a real-world data set to validate the efficiency of solving the MISOCP reformulation. 

Our computational experiments demonstrate that accounting for congestion costs and leads to the preservation of larger hub capacity levels, and the resulting hub network topology tends to exhibit a reduced number of dispersed hubs. We also observe a trade-off between decreasing transportation costs on the hub-hub link and increasing congestion costs due to consolidation.

\bibliographystyle{unsrtnat}
\bibliography{references}  

\begin{thebibliography}{16}
\providecommand{\natexlab}[1]{#1}
\providecommand{\url}[1]{\texttt{#1}}
\expandafter\ifx\csname urlstyle\endcsname\relax
  \providecommand{\doi}[1]{doi: #1}\else
  \providecommand{\doi}{doi: \begingroup \urlstyle{rm}\Url}\fi

\bibitem[Shen et~al.(2021)Shen, Liang, and Shen]{shen2021reliable}
Hao Shen, Yong Liang, and Zuo-Jun~Max Shen.
\newblock Reliable hub location model for air transportation networks under
  random disruptions.
\newblock \emph{Manufacturing \& Service Operations Management}, 23\penalty0
  (2):\penalty0 388--406, 2021.

\bibitem[Ernst and Krishnamoorthy(1999)]{ernst1999solution}
Andreas~T Ernst and Mohan Krishnamoorthy.
\newblock Solution algorithms for the capacitated single allocation hub
  location problem.
\newblock \emph{Annals of operations Research}, 86\penalty0 (0):\penalty0
  141--159, 1999.

\bibitem[Yaman et~al.(2012)Yaman, Karasan, and Kara]{yaman2012release}
Hande Yaman, Oya~Ekin Karasan, and Bahar~Y Kara.
\newblock Release time scheduling and hub location for next-day delivery.
\newblock \emph{Operations research}, 60\penalty0 (4):\penalty0 906--917, 2012.

\bibitem[Wu et~al.(2023)Wu, Herszterg, Savelsbergh, and Huang]{wu2023service}
Haotian Wu, Ian Herszterg, Martin Savelsbergh, and Yixiao Huang.
\newblock Service network design for same-day delivery with hub capacity
  constraints.
\newblock \emph{Transportation Science}, 57\penalty0 (1):\penalty0 273--287,
  2023.

\bibitem[Klincewicz(1998)]{klincewicz1998hub}
John~G Klincewicz.
\newblock Hub location in backbone/tributary network design: a review.
\newblock \emph{Location Science}, 6\penalty0 (1-4):\penalty0 307--335, 1998.

\bibitem[Khaniyev et~al.(2020)Khaniyev, Elhedhli, and
  Erenay]{khaniyev2020spatial}
Taghi Khaniyev, Samir Elhedhli, and Fatih~Safa Erenay.
\newblock Spatial separability in hub location problems with an application to
  brain connectivity networks.
\newblock \emph{INFORMS Journal on Optimization}, 2\penalty0 (4):\penalty0
  320--346, 2020.

\bibitem[Campbell and O'Kelly(2012)]{campbell2012twenty}
James~F Campbell and Morton~E O'Kelly.
\newblock Twenty-five years of hub location research.
\newblock \emph{Transportation Science}, 46\penalty0 (2):\penalty0 153--169,
  2012.

\bibitem[O’Kelly and Bryan(1998)]{o1998hub}
Morton~E O’Kelly and DL~Bryan.
\newblock Hub location with flow economies of scale.
\newblock \emph{Transportation research part B: Methodological}, 32\penalty0
  (8):\penalty0 605--616, 1998.

\bibitem[Alumur et~al.(2018)Alumur, Nickel, Rohrbeck, and Saldanha-da
  Gama]{alumur2018modeling}
Sibel~A Alumur, Stefan Nickel, Brita Rohrbeck, and Francisco Saldanha-da Gama.
\newblock Modeling congestion and service time in hub location problems.
\newblock \emph{Applied Mathematical Modelling}, 55:\penalty0 13--32, 2018.

\bibitem[Elhedhli and Wu(2010)]{elhedhli2010lagrangean}
Samir Elhedhli and Huyu Wu.
\newblock A lagrangean heuristic for hub-and-spoke system design with capacity
  selection and congestion.
\newblock \emph{INFORMS Journal on Computing}, 22\penalty0 (2):\penalty0
  282--296, 2010.

\bibitem[Bayram et~al.(2023)Bayram, Y{\i}ld{\i}z, and Farham]{bayram2023hub}
Vedat Bayram, Bar{\i}{\c{s}} Y{\i}ld{\i}z, and M~Saleh Farham.
\newblock Hub network design problem with capacity, congestion, and stochastic
  demand considerations.
\newblock \emph{Transportation Science}, 57\penalty0 (5):\penalty0 1276--1295,
  2023.

\bibitem[Kleinrock(2007)]{kleinrock2007communication}
Leonard Kleinrock.
\newblock \emph{Communication nets: Stochastic message flow and delay}.
\newblock Courier Corporation, 2007.

\bibitem[Campbell et~al.(2005)Campbell, Ernst, and
  Krishnamoorthy]{campbell2005hub}
James~F Campbell, Andreas~T Ernst, and Mohan Krishnamoorthy.
\newblock Hub arc location problems: part ii—formulations and optimal
  algorithms.
\newblock \emph{Management Science}, 51\penalty0 (10):\penalty0 1556--1571,
  2005.

\bibitem[Rostami et~al.(2022)Rostami, Chitsaz, Arslan, Laporte, and
  Lodi]{rostami2022single}
Borzou Rostami, Masoud Chitsaz, Okan Arslan, Gilbert Laporte, and Andrea Lodi.
\newblock Single allocation hub location with heterogeneous economies of scale.
\newblock \emph{Operations Research}, 70\penalty0 (2):\penalty0 766--785, 2022.

\bibitem[Alizadeh and Goldfarb(2003)]{alizadeh2003second}
Farid Alizadeh and Donald Goldfarb.
\newblock Second-order cone programming.
\newblock \emph{Mathematical programming}, 95\penalty0 (1):\penalty0 3--51,
  2003.

\bibitem[G{\"u}nl{\"u}k and Linderoth(2008)]{gunluk2008perspective}
Oktay G{\"u}nl{\"u}k and Jeff Linderoth.
\newblock Perspective relaxation of mixed integer nonlinear programs with
  indicator variables.
\newblock In \emph{International Conference on Integer Programming and
  Combinatorial Optimization}, pages 1--16. Springer, 2008.

\end{thebibliography}






\end{document}